\documentclass[12pt,a4paper]{article}
\usepackage[english]{babel}

\usepackage{latexsym}
\usepackage{amsfonts}
\usepackage{amssymb}
\usepackage{calrsfs}
\usepackage{amsmath, amsthm}

\newtheorem{Cor}{Corollary}
\newtheorem{note}{Remark}
\newtheorem{prop}{Proposition}
\newtheorem{lemma}{Lemma}

\newcommand{\N}{\mathbb N}

\begin{document}

\footnotetext {This work was completed with the support of  Russian Foundation for Basic Research (grant 14-01-00349-a, grant 13-01-12476-OFI-M2, grant 13-01-00384-a)}

\begin{center}
{\bf Examples of very unstable linear partial functional differential equations}

R. S. Ismagilov, N. A. Rautian, V. V. Vlasov

\end{center}

\begin{abstract}
We consider the examples of partial functional differential equations with delay in the Laplacian. First of these equations is linear parabolic equation, the second one is linear hyperbolic equation, third equation is perturbed hyperbolic equation with delay. We show that there are the sequence of eigenvalues in both cases with real parts tends to plus infinity.
\end{abstract}

\begin{center}
{\bf 1. Introduction}
\end{center}

Nowdays there exists many works devoted to the researching and comparison of different models of diffusions and heat conductions in media with memory (\cite{1}, \cite{7}). Most of these models use Maxwell-Cattaneo hyperbolic regularization of heat equation (Maxwell-Cattaneo equation) 
\begin{equation}\label{MW}
\frac{{\partial T(x,t)}}{{\partial t}} + \tau \frac{{\partial {T^2}(x,t)}}{{\partial {t^2}}} = \lambda  \cdot \Delta T(x,t).
\end{equation}
The left part of this equation is the first order Taylor expansion of the so called time delayed heat conduction equation
\begin{equation}\label{DE}
\frac{{\partial T(x,t + \tau )}}{{\partial t}} = \lambda  \cdot \Delta T(x,t).
\end{equation}
where $\tau  > 0$  corresponds to the time delay between cause and effect. The fact that time delayed heat conduction equation is more reasonable model was first noticed by Maxwell \cite{9}. There was  many investigations and numerical examples for a long time period showed that the Maxwell-Cattaneo equation (1) is a good approximation of the equation (2). But in the last decade of the twenty century and in the first decade of present century investigations in the field of non-stationary heat transfer processes became really widespread. The reason for it is a development of technologies.  There were many new models of heat conductivity, among which the ballistic-diffusive heat conduction model \cite{7}, model of hyperbolic self-consistent problem of heat transfer in rapid solidification of supercooled liquid, model of heat propagation dynamics in thin silicon layers etc. But all these models based on Maxwell-Cattaneo equation.

Our main purpose is to demonstrate that there exists qualitative difference between spectra of the time delayed equations and its hyperbolic regularizations and corresponding equations without delay. Namely, we will present the examples of the partial delay equations which  spectra have the sequence of  eigenvalues  ${\lambda _n}$  such that $\operatorname{Re} {\lambda _n} \to  + \infty $. We will call such equations unstable. In turn the spectra of the symbols of hyperbolic and parabolic equation lies in the left part $\left\{ {\lambda :\operatorname{Re} \lambda  < \omega }, \omega \in \mathbb R_+ \right\}$ of complex plane, thus hyperbolic and parabolic equations are stable in the sense defined above. Thus we will show that the spectra of the symbols of hyperbolic and parabolic equations seriously different from the spectra of partial functional differential equations.

We show motivated by the simple looking linear parabolic and hyperbolic equations with delay in Laplacian operator that initial value problems for these equations are awfully unstable. The heat equation with delay was considered earlier in [1]. It was shown in [1] that initial problem for this equation can be solved in the carefully chosen  Frechet space. Moreover it was shown in [1] that there exists a sequence of eigenvalues ${\lambda _n} = {x_n} + i{y_n}$ such that ${x_n} \to  + 0$ ($n \to  + \infty $). Thus the authors obtained the lack of exponential dichotomy. They note that heat equation with delay arises when we consider random movement of a biological species and when we assume spatial movement of the species is delayed.

\begin{center}
{\bf 2. Examples}
\end{center}

{\bf Example 1.}
We consider the heat equation with delay of the following form:
\begin{equation}\label{eq:1}
{u_{t}} = {u_{xx}}(t - h,x), \quad t>0, \, 0 < x < \pi, \, h>0 
\end{equation}
with Dirihlet boundary conditions
\begin{equation}\label{eq:2}
{\left. u \right|_{x = 0}} = {\left. u \right|_{x = \pi }} = 0. 
\end{equation}
The main purpose of our considerations is to study the spectrum distribution of the symbol of equation \eqref{eq:1}. In oder to do this we will look for the solution of the equation \eqref{eq:1} in the form
$$
u(t,x) = \sum\limits_{n = 1}^\infty  {{T_n}(t)\sin nx}, 
$$
using the Fourier method.
Then we obtain infinite number of ordinary delay equations 
\begin{equation}\label{eq:3}
{T'_n}(t) =  - {n^2}{T_n}(t - h), \quad n\in \mathbb N,
\end{equation}
from the equation \eqref{eq:1}.The following equations 
\begin{equation}\label{eq:9}
\lambda  + {n^2}{e^{ - \lambda h}}=0, \quad n\in \N, \quad \lambda= x+iy,
\end{equation}
are the symbols (characteristic quasipolinomials) of the equations \eqref{eq:3}. Let us  put $h=1$ for the symplicity.

{\bf Example 2.}
Now let us consider wave equation with delay
\begin{equation}\label{eq:5}
{u_{tt}} = {u_{xx}}(t - h,x), \quad t>0, \, 0 < x < \pi, \, h>0,
\end{equation}
with Dirihlet boundary conditions
\begin{equation}\label{eq:6}
{\left. u \right|_{x = 0}} = {\left. u \right|_{x = \pi }} = 0. 
\end{equation}
We will look for the solution of equation \eqref{eq:5} in the following form 
$$
u(t,x) = \sum\limits_{n = 1}^\infty  {{T_n}(t)\sin nx},
$$
using the Fourier method. Thus we obtain the infinite number of ordinary delay equations 
\begin{equation}\label{eq:7}
{T''_n}(t) =  - {n^2}{T_n}(t - h), \quad n\in \mathbb N.
\end{equation}
The following quasipolinomials 
\begin{equation}\label{eq:8}
\lambda^2  + {n^2}{e^{ - \lambda h}} = 0, \quad n\in \N, \quad \lambda= x+iy
\end{equation}
are the symbols of the equations \eqref{eq:8}. Let us  put $h=1$ for the symplicity.

{\bf Example 3.}
Now let us consider perturbed wave equation with delay
\begin{equation}\label{eq:55}
{u_{tt}} = {u_{xx}}(t,x)+{u_{xx}}(t - h,x), \quad t>0, \, 0 < x < \pi, \, h>0,
\end{equation}
with Dirihlet boundary conditions
\begin{equation}\label{eq:66}
{\left. u \right|_{x = 0}} = {\left. u \right|_{x = \pi }} = 0. 
\end{equation}
We will look for the solution of equation \eqref{eq:5} in the following form 
$$
u(t,x) = \sum\limits_{n = 1}^\infty  {{T_n}(t)\sin nx},
$$
using the Fourier method. Thus we obtain the infinite number of ordinary delay equations 
\begin{equation}\label{eq:77}
{T''_n}(t) =  - {n^2}\left({T_n}(t)+{T_n}(t - h)\right), \quad n\in \mathbb N.
\end{equation}
The following quasipolinomials 
\begin{equation}\label{eq:88}
\lambda^2  + {n^2}(1+{e^{ - \lambda h}}) = 0, \quad n\in \N, \quad \lambda= x+iy
\end{equation}
are the symbols of the equations \eqref{eq:8}. Let us put $h=1$ for the symplicity.

\begin{center}
{\bf 2. Statements of the results and Proofs.}
\end{center}

First we shall prove that the equations \eqref{eq:9}, \eqref{eq:8}, \eqref{eq:88} are very unstable. Equations \eqref{eq:9}, \eqref{eq:8} can be written in the following form 
\begin{equation}\label{eq:111}
\lambda  + b \cdot \ln \lambda  - w = 0,
\end{equation}
where $w = 2\ln n + i\pi $ and constant $b = 1$ for the equation \eqref{eq:9} and $b = 2$ for the equation \eqref{eq:8}. Here we consider such brunch of that $\ln \lambda  = \ln |\lambda | + i \cdot \arg \lambda $, $\arg \lambda  \in ( - \pi ,\pi )$.

\begin{lemma}\label{1}
Let us consider $w \in \mathbb{C}$, $\operatorname{Re} w > 0$, $r = |w|$. If $r$ is sufficiently large then there exists the unique solution $\lambda  = \lambda (w)$ of the equations \eqref{eq:9}, \eqref{eq:8} in the circle $|\lambda  - w| < r/2$.
\end{lemma}

{\bf Proof. }
Let us consider the equation \eqref{eq:111} and equation $\lambda  - w = 0$. We have the inequality $|\lambda | \leqslant 3/2 \cdot r$ on the circle $|\lambda  - w| = r/2$ and hence we obtain inequality $|b\ln \lambda | \leqslant b\ln r + C$ where constant $C$ depends on $b$. So we have $|b\ln \lambda | < |\lambda  - w|$ on the circle $|\lambda  - w| = r/2$ and due to Rouche theorem the equation \eqref{eq:9} and \eqref{eq:8} has the unique solution $\lambda  = \lambda (w)$ in the circle $|\lambda  - w| < r/2$.

\begin{Cor}
We have $\operatorname{Re} \lambda  \to  + \infty $ when $\operatorname{Re} w \to  + \infty $.
\end{Cor}

Due to corollary we obtain that equations \eqref{eq:9}, \eqref{eq:8} are very unstable in the following sense: there exists such solution ${\lambda _n} = {x_n} + i{y_n}$ of the equation \eqref{eq:9} or  \eqref{eq:8}  that ${x_n} \to  + \infty$ ($n \to  + \infty $).

\begin{note}
If we substitute $n^2$ by $n^{\theta}$ $(\theta>0)$ in the equations  \eqref{eq:9}, \eqref{eq:8} the results of Lemma \ref{1} and Corollary \ref{1} will be valid. Due to this fact the equations 
$$
\begin{gathered}
  {u_t} = \Delta u(t - h,x), \hfill \\
  {u_{tt}} = \Delta u(t - h,x),\quad x \in G \subset {\mathbb{R}^N}, \quad t>0,\quad h>0  \hfill \\ 
\end{gathered} 
$$
with Dirithlet boundary conditions
$$
{\left. u \right|_{\partial G}} = 0,
$$
where $G$ is bounded domain with smooth boundary are also unstable.
\end{note}

Moreover if we substitute Laplacian in these equations by more general elliptic selfajoint operator of oder $2m$ in bounded domain with smooth boundary the corresponding equations will be also unstable.

In the following propositions 1 and 2 we present the concrete sequences of zeroes $\lambda_n= x_n+iy_n$  of the equations \eqref{eq:9} and \eqref{eq:8} such that $x_n\to +\infty$.

\begin{prop}
There is a family of solutions $\lambda_n= x_n+iy_n$ of the equation \eqref{eq:9} such that ${x_n} \sim \ln {n^2} - \ln \ln {n^2}$, for  $n\to +\infty$.
\end{prop}
{\bf Proof. } We can write  \eqref{eq:9} like the following system extracting real and imaginary parts
\begin{equation}\label{eq:10}
\left\{ \begin{gathered}
  {e^x}\left( {x\cos y - y\sin y} \right) =  - {n^2}, \hfill \\
  y\cos y + x\sin y = 0. \hfill \\ 
\end{gathered}  \right. 
\end{equation}
If $y=0$ then \eqref{eq:10} has unique sulution $(0,0)$ for $n=0$. Let us put $y \ne 0$ then we have
$$
x =  - \frac{{y\cos y}}{{\sin y}}
$$
and 
$$
{e^{ - \frac{{y\cos y}}{{\sin y}}}}\frac{y}{{\sin y}} = {n^2}. 
$$
Note that $x =g(y)=  - \frac{{y\cos y}}{{\sin y}}\to +\infty$ for $y\to \pi-0$. If we put $y = \pi  - \delta$, \, $\delta>0$ the equation (3) has the following form
$$
{e^{\frac{{(\pi  - \delta )\cos \delta }}{{\sin \delta }}}}\frac{{(\pi  - \delta )}}{{\sin \delta }} = {n^2},
$$
that equivalent to the equation
$$
{e^{\left( {\frac{\pi }{\delta } - 1} \right)}}\left( {\frac{\pi }{\delta } - 1} \right) = {n^2}
$$
as $\delta\to 0$. Denote by $\theta  = \displaystyle\frac{\pi}{{\delta }}-1$ then (4) has the form
$$
\theta {e^\theta } = {n^2}.
$$
Using the results from the monograph of M.V. Fedoryuk ([2], pp. 51--52 in Russian) we obtain the following asymptotic representation for $\theta$:
$$
\theta  = {\ln {n^2}} - {\ln \ln {n^2}}+O\left( {\frac{{\ln \ln {n^2}}}{{\ln {n^2}}}} \right), \quad n\to +\infty.
$$
Thus we have 
$$
\frac{\pi }{{\delta }} = {\ln {n^2}} - {\ln \ln {n^2}}+O\left( {\frac{{\ln \ln {n^2}}}{{\ln {n^2}}}} \right), \quad n\to +\infty
$$
and
$$
\delta  \sim \sin \delta  \sim \frac{\pi }{{\ln {n^2}} - {\ln \ln {n^2}}}, \quad n\to +\infty.
$$
Then
$$
{y_n} \sim \pi \left( {1 - \frac{1}{{{\ln {n^2}} - {\ln \ln {n^2}}}}} \right), \quad n\to +\infty
$$
and 
\begin{multline*}
{x_n} =  - {y_n}\frac{{\cos {y_n}}}{{\sin {y_n}}} \sim  - \pi \left( {1 - \frac{1}{\ln {n^2} - \ln \ln {n^2}}} \right)\frac{{\cos \pi \left( {1 - \frac{1}{{{\ln {n^2}} - {\ln \ln {n^2}}}}} \right)}}{{\sin \pi \left( {1 - \frac{1}{{{\ln {n^2}} - {\ln \ln {n^2}}}}} \right)}} \sim \\ \sim {\ln {n^2}} - {\ln \ln {n^2}}\to  + \infty, \quad n\to +\infty.
\end{multline*}

\begin{prop}
There exists a family of solutions $\lambda_n= x_n+iy_n$ of the equation \eqref{eq:8} such that ${x_n} \sim 2(\ln (n/2) - \ln \ln (n/2))$, for  $n\to +\infty$.
\end{prop}
{\bf Proof. } Let us extract the real and imaginary parts of the equation \eqref{eq:8}:
\begin{equation}\label{eq:11}
\left\{ \begin{gathered}
  {x^2} - {y^2} + {n^2}{e^{ - x}}\cos y = 0 \hfill \\
  2xy - {n^2}{e^{ - x}}\sin y = 0 \hfill \\ 
\end{gathered}  \right.
\end{equation}
We deduce from \eqref{eq:11} the following system
\begin{equation}\label{eq:12}
\left\{ \begin{gathered}
  {\text{ctg}}y = \frac{{{y^2} - {x^2}}}{{2xy}} = \frac{1}{2}\left( {\frac{y}{x} - \frac{x}{y}} \right), \hfill \\
  \left( {{x^2} + {y^2}} \right){e^x} = {n^2}. \hfill \\ 
\end{gathered}  \right.
\end{equation}
We put $t = \displaystyle\frac{x}{y}$. Then we have for $t$ the following equation
$$
2{\text{ctg}}y = \frac{1}{t} - t,
$$
hence
$$
{t^2} + 2{\text{ctg}}y \cdot t - 1 = 0,
$$
$$
t =  - {\text{ctg}}y \pm \sqrt {{\text{ct}}{{\text{g}}^2}y + 1}  = \frac{{ - \cos y \pm 1}}{{\sin y}}.
$$
Thus we have
$$
x = y\left( {\frac{{ - \cos y \pm 1}}{{\sin y}}} \right).
$$
Denote by $y = \pi  - \delta $ and consider small enough $\delta>0$. Then we obtain
\begin{equation}\label{eq:112}
x = (\pi  - \delta )\frac{{1 + \cos \delta }}{{\sin \delta }} =\frac{{2\pi }}{\delta } - 2 - \frac{{\pi \delta }}{6} + \frac{{{\delta ^2}}}{6} + ({\delta ^2}), \quad \delta  \to +0. 
\end{equation}
Then we obtain from the second equation of the system \eqref{eq:12} the following equation
\begin{multline}\label{eq:13}
{(\pi  - \delta )^2}\left[ {1 + {{\left( {\frac{{1 - \cos (\pi  - \delta )}}{{\sin \delta }}} \right)}^2}} \right]{e^{(\pi  - \delta )\left( {\frac{{1 - \cos (\pi  - \delta )}}{{\sin \delta }}} \right)}} \approx \left( {\frac{{4{\pi ^2}}}{{{\delta ^2}}} - \frac{{8\pi }}{\delta } + 4} \right){e^{\frac{{2\pi }}{\delta } - 2}}= {n^2}, \\ \delta  \to +0.
\end{multline}
Denote $\theta  = \displaystyle\frac{{2\pi }}{\delta }$, $\theta\to  + \infty $, $\left( {\delta  \to  + 0} \right)$ then the we can write the equation \eqref{eq:13} in the following form
$$
\left( {{\theta}-2} \right)^2e^{\theta -2 } = {n^2}
$$
Denote $\eta  = \theta  - 2$, then we obtain the following equation
$$
{\eta ^2}{e^\eta } = {n^2}.
$$
Using the results from the monograph of M.V. Fedoryuk ([1], pp. 51--52 in Russian) we have
$$
\eta  = 2\left( {\ln {\frac{n}{2}} - \ln \ln{\frac{n}{2}} } \right) + O\left( {\frac{{\ln \ln{\frac{n}{2}}}}{{\ln {\frac{n}{2}}}}} \right), \quad n\to +\infty
$$
Hence we obtain the following asyptotic representations:
$$
\frac{{2\pi }}{\delta } =\theta = 2\left( {\ln {\frac{n}{2}} - \ln \ln{\frac{n}{2}} } \right) + O\left( {\frac{{\ln \ln{\frac{n}{2}}}}{{\ln {\frac{n}{2}}}}} \right)+2, \quad \delta  \to 0 +, \, n\to +\infty.
$$
$$
\delta  = \frac{2\pi }{{\left( {\ln \displaystyle\frac{n}{2} - \ln \ln \frac{n}{2}} \right) + 2 + O\left( {\displaystyle\frac{{\ln \ln {\frac{n}{2}}}}{{\ln {\frac{n}{2}}}}} \right)}}, \quad n\to +\infty,
$$
Thus we obtain the following asymptotic representations from the representation \eqref{eq:112} :
$$
{x_n} =2\left( {\ln {\frac{n}{2}} - \ln \ln{\frac{n}{2}} } \right) + O\left( {\frac{{\ln \ln{\frac{n}{2}}}}{{\ln {\frac{n}{2}}}}} \right), \quad n\to +\infty.
$$

Consider the equations \eqref{eq:88} for ${\lambda ^2} + {n^2}\left( {1 + {e^{ - \lambda }}} \right) = 0,$ $n \in \mathbb{N}$.

\begin{lemma}
There exists a sequence  $ {\lambda _n}=x_n+iy_n$ of the solutions of the equations \eqref{eq:88} such that $\operatorname{Re}{\lambda _n} =x_n\to +\infty$ for $n\to +\infty$. Moreover 
$$
x_n ={\ln {n} - \ln \ln{n} } + O\displaystyle\left( {\frac{{\ln \ln{n}}}{{\ln {n}}}} \right), \quad n \to \infty.
$$
\end{lemma}

{\bf Proof. }
We need the equation $x{e^x}=t$  where $t>0$, $x>0$. This equation has a unique solution $x = \Phi (t)$.  Hence   $e^x = \displaystyle\frac{t}{{\Phi (t)}}$  and  $\Phi (t) ={\ln {t} - \ln \ln{t} } + O\displaystyle\left( {\frac{{\ln \ln{t}}}{{\ln {t}}}} \right)\simeq \ln t$ ($t \to \infty $). We put $\lambda  = x+iy$. Then the equation \eqref{eq:55} may be written like a system 
\begin{equation}\label{eq:112*}
\left\{ \begin{gathered}
  {x^2} - {y^2} + {n^2}(1 + {e^{ - x}}\cos y) = 0, \hfill \\
  2xy - {n^2}{e^{ - x}}\sin y = 0. \hfill \\ 
\end{gathered}  \right.
\end{equation}
The second equation in \eqref{eq:112*} we shall rewrite in the following way:
$$
x{e^x} = t, \quad t = \frac{{{n^2}\sin y}}{{2y}}
$$
Hence
$$
x = \Phi (t), \quad {e^x} = \frac{t}{{\Phi (t)}}.
$$
Substituting it into the first equation of the system \eqref{eq:112*} we obtain the equation for function $\Phi (t)$: ${\Phi ^2}(t) - {y^2} + {n^2} + {n^2}\Phi (t)\cos y/t = 0$. Solving this equation we obtain 
$$
\Phi (t) = \Phi \left(\frac{{{n^2}\sin y}}{{2y}}\right) =  - \frac{{y\cos y - \sqrt {{y^2} - {n^2}{{\sin }^2}y} }}{{\sin y}}
$$
Let us denote
$$
{U_n}(y) = \Phi \left(\frac{{{n^2}\sin y}}{{2y}}\right) 
+ \frac{{y\cos y - \sqrt {{y^2} - {n^2}{{\sin }^2}y} }}{{\sin y}}.
$$
Consider the integers $n$ such that $\cos n > \alpha  > 0$  and  $\cos (n + 1) <  - 1/4$. It is possible to show that
There is infinite number of such integers. For these numbers we have   $\sin n > \alpha$, $\sin (n + 1) > \alpha $
where  $\alpha >0$. Consider function ${U_n}(y)$ for  $y \in \left[ {n,n + 1} \right]$. We will show that numbers $U(n)$  and $U(n + 1)$
have different signs. We have   $U(n) = \Phi \left( {\displaystyle\frac{{n\sin n}}{2}} \right) > 0$. From the other hand we have 
$$
{U_n}(n + 1) < \Phi \left( {\frac{{{{(n + 1)}^2}\sin (n + 1)}}{{2n + 2}}} \right) + \frac{{(n + 1)\cos (n + 1)}}{{\sin (n + 1)}} < {c_1}\ln n - {c_2}n, \quad {c_1}, {c_2} > 0.
$$

Here we used the inequality  $\Phi (t) < t$ and also the inequalities $\cos (n + 1) < -1/4$, $\sin (n + 1) > 0$.
So we obtain that  ${U_n}(n + 1) < 0$ if $n$  is sufficiently large. Hence, the equation ${U_n}(y) = 0$ has the solution  ${y_n} \in (n,n + 1)$. Using the inequality $\sin {y_n} > \alpha $ we obtain that  ${x_n} = \Phi \left( {\displaystyle\frac{{{n^2}\sin {y_n}}}{{2{y_n}}}} \right) \to +\infty $ and  $x_n ={\ln {n} - \ln \ln{n} } + O\displaystyle\left( {\frac{{\ln \ln{n}}}{{\ln {n}}}} \right)$ ($n \to \infty $).

\medskip
\noindent {\bf Remark.}
In comparison with hyperbolic case (equation \eqref{eq:88}), parabolic equation with delay 
\begin{equation}\label{eq:R}
{u_t} = {u_{xx}}(t,x) + {u_{xx}}(t - h,x), \quad 0<x<\pi,\quad t>0, \quad h>0
\end{equation}

is stable in the following sence: the semiplate $\left\{ {\lambda :\operatorname{Re} \lambda  > \omega } \right\}$  is free for any $\omega>0$ of eigenvalues $\lambda_n$. That is for arbitrary zeroes $\lambda_n=x_n+iy_n$ of characteristic quasipolinomials  
$$
\lambda  + {n^2}(1 + {e^{ - \lambda h}}) = 0, \quad n \in \mathbb{N} \quad 
$$
the real parts ${x_n} \leqslant \omega $.

It is relevant to note that equation \eqref{eq:R} can be written in abstract form 
\begin{equation}\label{eq:RA}
\frac{{du}}{{dt}} + {A^2}u(t) + {A^2}u(t - h) = 0 \quad h>0,
\end{equation}
where ${A^2}y =  - y''(x)$, $y(0)=y(\pi)=0$.

The equation \eqref{eq:RA} is the simplest case of the equations which were considered in many articles. We restrict ourselves and cite only articles [3]-[5]. The abstract parabolic functional differential equations with unbounded operator coefficients were considered in the articles [3]-[4].  The main part of these equations is the abstract parabolic equation
$$
\frac{{du}}{{dt}} + {A^2}u(t)=0.
$$
where $A^2$ is selfadjoint positive operator, having compact inverse. The correct solvability of functional differential equations mentioned above was obtained in weighted Sobolev spaces $W_{2,\gamma }^1\left( {{\mathbb{R}_ + },{A^2}} \right)$. Moreover it was shown in autonomous case in [3] (see lemma 2, proposition 3 and lemma 3 for details) that symbol of this equation (analogue of characteristic quasipolinomial) is invertible in the semiplate $\left\{ {\lambda :\operatorname{Re} \lambda  > \gamma } \right\}$. So in this situation there are no sequences of eigenvalues  $\lambda_n=x_n+iy_n$ such that $x_n\to +\infty$.

It is relevant to note that abstract functional differential equations having main part is abstract hyperbolic equation 
$$
\frac{{d^2u}}{{dt^2}} + {A^2}u(t)=0
$$
was considered in \cite{6}.

\begin{center}
{\bf 2. Concluding remarks.}
\end{center}

The simplest equation \eqref{eq:1}, \eqref{eq:5}, \eqref{eq:55} considered in this article can be written in the following abstract form:
\begin{equation}\label{eq:A1}
\frac{{du}}{{dt}} + {A^2}u(t - h) = 0,
\end{equation}
\begin{equation}\label{eq:A2}
\frac{{{d^2}u}}{{d{t^2}}} + {A^2}u(t - h) = 0,
\end{equation}
\begin{equation}\label{eq:A3}
\frac{{{d^2}u}}{{d{t^2}}} + {A^2}u(t) + {A^2}u(t - h) = 0.
\end{equation}
where $A^2$ is selfadjoint positive operator in the Hilbert space $H \equiv {L_2}\left( {0,\pi } \right)$ having compact inverse, ${A^2}y =  - y''(x)$, $y(0)=y(\pi)=0$. The examples 1-3 show that classical initial problems for these equations can't be solved in weighted Sobolev spaces $W_{2,\gamma }^n\left( {{\mathbb{R}_ + },{A^n}} \right)$. The understandable reason of this fact is the Laplace transforms of the functions from the space $W_{2,\gamma }^n\left( {{\mathbb{R}_ + },{A^n}} \right)$ is analytic in the semiplate  $\left\{ {\lambda :\operatorname{Re} \lambda  > \gamma } \right\}$.

Let us consider the spectrum of the equation (1) in one dimensional case under the following assumptions 
$$
\lambda  = 1,\quad \tau  = 1,\quad \Delta T =  T_{xx}(t,x),\quad T(t,0) = T(t, \pi ) = 0.
$$
Using the Fourier method
$$
T(t,x) = \sum\limits_{n = 1}^\infty  {{T_n}(t)\sin nx,}
$$
we obtain the following ordinary differential equations
$$
{T''_n}(t) + {T'_n}(t) =  - {n^2}{T_n}(t), \quad n \in \mathbb{N}.
$$
The characteristic polynomials have the form
$$
{\lambda ^2} + \lambda  + {n^2} = 0, \quad n \in \mathbb{N}
$$
Hence we have
$$
\lambda _n^ \pm  = \frac{{ - 1 \pm \sqrt {1 - 4{n^2}} }}{2} =  - \frac{1}{2} \pm in\left( {1 - \frac{1}{{8{n^2}}} + o\left( {\frac{1}{{{n^3}}}} \right)} \right)
$$
and the spectra of this problem is
$$
\sigma  = \bigcup\limits_{n = 1}^\infty  {{\lambda ^ \pm }_n}.
$$
At the same time the spectra of the problem (6), (7) which  coincides with the spectra of the equation (2) for $\lambda  = 1$, $\tau  = 1$ and it can be represented in the following way
\begin{equation}\label{eq:sigma}
\Sigma  = \overline {\bigcup\limits_{n = 1}^\infty  {\bigcup\limits_{k \in \mathbb{Z}}^{} {{\lambda _{nk}}} } }
\end{equation}
where ${\lambda _{nk}}$ are the zeroes of quasipolinomials (6) having the following asymptotic representations (for fixed n)
\begin{equation}\label{eq:lambda}
{\lambda _{nk}}^ \pm  = \ln ({n^2}) - \ln \left| {\frac{\pi }{2} + 2\pi k} \right| + \frac{\pi }{{2\pi k}} + o\left( {\frac{1}{k}} \right) \pm i\left[ {\frac{\pi }{2} + 2\pi k + O\left( {\frac{{\ln k}}{k}} \right)} \right].
\end{equation}
(see, for example monograph \cite{8}, chapter 4).

Formulas \eqref{eq:sigma} and \eqref{eq:lambda} shows that the structure of the spectra of the equation \eqref{eq:1} and \eqref{eq:2} are seriously differs. Example 1 wchich was considered in this paper confirm it. So it is naive to expect that properties of solutions of the equation \eqref{eq:1} and \eqref{eq:2} will be similar.

\vspace{0.5ex}

\begin{flushleft}

{\footnotesize R. S. Ismagilov

Department of Higher Mathematics

Bauman Moscow State Technical University 

ul. Baumanskaya 2-ya, 5/1, Moscow, 105005 

Russia}

{\footnotesize { \it Email address:} ismagil@serv.bmstu.ru}

\end{flushleft}

\vspace{0.5ex}

\begin{flushleft}

{\footnotesize N. A. Rautian

Department of Mechanics and Mathematics

Moscow Lomonosov State University

Vorobievi Gori, Moscow, 117234 

Russia}

{\footnotesize { \it Email address:} nrautian@mail.ru}

\end{flushleft}

\vspace{0.5ex}

\begin{flushleft}

{\footnotesize V. V. Vlasov

Department of Mechanics and Mathematics

Moscow Lomonosov State University

Vorobievi Gori, Moscow, 117234 

Russia}

{\footnotesize { \it Email address:} : vlasovvv@mech.math.msu.su}

\end{flushleft}


\begin{thebibliography}{1}

\bibitem{1}
{Rodrigues H.~M., Ou C., Wu J.}
\textit {A Partial Differential Equation with Delayed Diffusion Dynamics of Continuous, Discrete and Impulsive Systems.} 
\textbf (to appear)

\bibitem{2}
{Fedorjuk M.~V.}
\textit{Asymptotics, Integrals and Series}{Moscow "`Nauka"',} 1987




\bibitem{3}
{Vlasov V.~V.} 
\textit{On the solvability and properties of solutions of functional-differential equations in Hilbert spaces.}  Sbornik Mathem.
\textbf{186:8} (1995), 67--92.


\bibitem{4}
{Vlasov V.~V.} 
\textit{On the solvability and estimates of solutions of functional-differential equations in Sobolev spaces.}  Proc. Steklova Math. Inst. 
\textbf{227} (1999), 109--121.


\bibitem{5}
{Di Blasio G., Kunisch K., Sinestari E.} 
\textit{$L^2$--regularity for parabolic partial integrodifferential equations with delays in the highest order derivatives.} J. Math. Anal. Appl.
\textbf{102} (1984), 38--57.


\bibitem{6}
{Vlasov V.~V.,  Shmatov K.~I.} 
\textit{On the solvability of delayed hyperbolic equations in Hilbert spaces.}  Proc. Steklov Math. Inst.
\textbf{243} (2003), 127--137.

\bibitem{8}
{Bellman R., Cooke K.~L.}
\textit {Differential-Difference Equation.} New York Academic Press London, 1963.


\bibitem{7}
{Xu~M., Hu~H.} 
\textit {A ballistic-diffusive heat conduction model extracted from Boltzman transport equation.} Proc. of the Royal Soc. A, January 2011.


\bibitem{9}
{Maxwell J.~C.}
\textit{On the dynamical theory of gases.} Phil. Trans. Royal Soc. 
\textbf{157} (1867), 49--88.






\end{thebibliography}
\end{document}